\documentclass[11pt]{article}
\usepackage{amsfonts}
\usepackage{amsfonts}
\usepackage{amsfonts}
\usepackage{mathrsfs}
\usepackage{bm}
\usepackage{amssymb,amsmath}
\usepackage{cite}
 \allowdisplaybreaks

\catcode`!=11
\let\!int\int \def\int{\displaystyle\!int}
\let\!lim\lim \def\lim{\displaystyle\!lim}
\let\!sum\sum \def\sum{\displaystyle\!sum}
\let\!sup\sup \def\sup{\displaystyle\!sup}
\let\!inf\inf \def\inf{\displaystyle\!inf}
\let\!cap\cap \def\cap{\displaystyle\!cap}
\let\!max\max \def\max{\displaystyle\!max}
\let\!min\min \def\min{\displaystyle\!min}
\let\!frac\frac \def\frac{\displaystyle\!frac}
\catcode`!=12

\let\oldsection\section
\renewcommand\section{\setcounter{equation}{0}\oldsection}

\allowdisplaybreaks
\def\pf{\it{Proof.}\rm\quad}

\def\N{\mathbb{N}}\def\Z{\mathbb{Z}}

\newtheorem{thm}{Theorem}[section]
\newtheorem{lem}[thm]{Lemma}
\newtheorem{cor}[thm]{Corollary}

\begin{document}
\title {\bf Cubic alternating harmonic number sums}
\author{
{Ce Xu\thanks{Corresponding author. Email: 15959259051@163.com (C. Xu)}}\\[1mm]
\small  School of Mathematical Sciences, Xiamen University\\
\small Xiamen
361005, P.R. China}

\date{}
\maketitle \noindent{\bf Abstract } A recent paper of A. Sofo proves some results about sums of products of quadratic alternating harmonic numbers and reciprocal binomial coefficients. In this paper, we extend his result to cubic alternating harmonic number sums and develop new closed form representations of sums of cubic alternating harmonic numbers and reciprocal binomial coefficients. Some interesting (known or new) illustrative special cases as well as immediate consequences of the main results are also considered.
\\[2mm]
\noindent{\bf Keywords} Harmonic numbers; Riemann zeta functions; Binomial coefficients; multiple harmonic (star) sums; polylogarithm functions.
\\[2mm]
\noindent{\bf AMS Subject Classifications (2010):} 05A10; 05A19; 11B65; 11M06; 11M32

\section{Introduction}
In a recent paper \cite{S2015}, A. Sofo prove some results on sums of products of alternating quadratic harmonic numbers and reciprocal binomial coefficients of the form
\begin{eqnarray}
{\overline W_k}\left( {1,1;p} \right): = \sum\limits_{n = 1}^\infty  {\frac{{H_n^2}}{{{n^p}\left( {\begin{array}{*{20}{c}}
   {n + k}  \\
   k  \\
\end{array}} \right)}}{{\left( { - 1} \right)}^{n + 1}}}\quad (k\in\N)\end{eqnarray}
for $p=0$ or $1$. The generalized $n$-th harmonic number of order $k$, $H^{(k)}_n$, is defined for positive integers $n$ and $r$ as (\cite{B1985,B1989,XMZ2016,X2017})
\begin{eqnarray}H^{(k)}_n:=\sum\limits_{j=1}^n\frac {1}{j^k}\quad (n, k \in \N)\end{eqnarray}
where the empty sum $H^{(m)}_0$ is conventionally understood to be zero, and $H_n:=H^{(1)}_n$.
In this paper we will develop identities, closed form representations of alternating
quadratic and cubic harmonic numbers and reciprocal binomial coefficients of the form:
\begin{eqnarray}
{{\overline W}_k}\left( {1,2;p} \right): = \sum\limits_{n = 1}^\infty  {\frac{{{H_n}H_n^{\left( 2 \right)}}}{{{n^p}\left( {\begin{array}{*{20}{c}}
   {n + k}  \\
   k  \\
\end{array}} \right)}}{{\left( { - 1} \right)}^{n + 1}}} \ {\rm and}\  {{\overline W}_k}\left( {{{\left\{ 1 \right\}}_3};p} \right): = \sum\limits_{n = 1}^\infty  {\frac{{H_n^3}}{{{n^p}\left( {\begin{array}{*{20}{c}}
   {n + k}  \\
   k  \\
\end{array}} \right)}}{{\left( { - 1} \right)}^{n + 1}}} ,
\end{eqnarray}
for $p=0$ and $1$. Here the notation $\{\}_p$ means that the sequence in the bracket is repeated $p$-times.
The generalized sums of products of alternating harmonic numbers and reciprocal binomial coefficients ${\overline W_k}\left( {{m_1},{m_2}, \cdots ,{m_r};p} \right)$ are defined by
\begin{eqnarray}{\overline W_k}\left( {{m_1},{m_2}, \cdots ,{m_r};p} \right): = \sum\limits_{n = 1}^\infty  {\frac{{H_n^{\left( {{m_1}} \right)}H_n^{\left( {{m_2}} \right)} \cdots H_n^{\left( {{m_r}} \right)}}}{{{n^p}\left( {\begin{array}{*{20}{c}}
   {n + k}  \\
   k  \\
\end{array}} \right)}}}(-1)^{n+1} ,\:\left( {p \in {\N_0},k,r,{m_i} \in \N} \right),\end{eqnarray}
where $\N:=\{1,2,\ldots\}\ {\rm and}\ \N_0:=\N\cup \{0\}=\{0,1,2,\ldots\}$.

While there are many results for sums of harmonic numbers with positive terms.  Many
harmonic number sums can be expressed in terms of a linear rational combination of classical
Riemann zeta values and harmonic numbers. For example we know that \cite{FS1998,X2016,X2017}
\[\sum\limits_{n = 1}^\infty  {\frac{{H_n^3}}{{{n^4}}}}  =
\frac{{231}}{{16}}\zeta (7) - \frac{{51}}{4}\zeta (3)\zeta (4) +
2\zeta (2)\zeta (5),\]
\[\sum\limits_{n = 1}^\infty  {\frac{{[H^{(2)} _n]^2}}{{{n^5}}}}  =  - \frac{{1069}}{{36}}\zeta \left( 9 \right) + \frac{4}{3}{\zeta ^3}\left( 3 \right) + 7\zeta \left( 2 \right)\zeta \left( 7 \right) - \frac{4}{3}\zeta \left( 3 \right)\zeta \left( 6 \right) + \frac{{33}}{2}\zeta \left( 4 \right)\zeta \left( 5 \right),\]
\[\sum\limits_{n = 1}^\infty  {\frac{{[H^{(2)}_n]^2H^{(3)}_n}}{{{n^2}}}}  =  - \frac{{617}}{{72}}\zeta \left( 9 \right) + {\zeta ^3}\left( 3 \right) + \frac{{91}}{8}\zeta \left( 2 \right)\zeta \left( 7 \right) - \frac{{17}}{4}\zeta \left( 4 \right)\zeta \left( 5 \right) - \frac{{329}}{{84}}\zeta \left( 3 \right)\zeta \left( 6 \right)\]
and from \cite{SX2017,XMZ2016}
\begin{align*}
\sum\limits_{n = 1}^\infty  {\frac{{{H_n}H_n^{(2)}}}{{n\left( {\begin{array}{*{20}{c}}
   {n + k}  \\
   k  \\
\end{array}} \right)}}}  = \sum\limits_{r = 1}^k {{{\left( { - 1} \right)}^{r + 1}}\left( {\begin{array}{*{20}{c}}
   k  \\
   r  \\
\end{array}} \right)} \left\{ {\begin{array}{*{20}{c}}
   {2\zeta \left( 4 \right) + 2\zeta \left( 3 \right){H_{r - 1}} + \frac{1}{2}\zeta \left( 2 \right)H_{r - 1}^2 + \sum\limits_{i = 1}^{r - 1} {\frac{{{H_i}}}{{{i^3}}}} }  \\
   { - \frac{1}{2}\zeta \left( 2 \right)H_{r - 1}^{(2)} - \frac{1}{2}\sum\limits_{i = 1}^{r - 1} {\frac{{H_i^2 + H_i^{(2)}}}{{{i^2}}}}  - \sum\limits_{i = 1}^{r - 1} {\frac{1}{i}\sum\limits_{j = 1}^i {\frac{{{H_j}}}{{{j^2}}}} } }  \\
\end{array}} \right\},
\end{align*}
there are fewer results for sums of the type (1.4). Here the classical Riemann zeta function defined by (\cite{B1985,B1989,SC2012})
\begin{eqnarray}\zeta \left( p \right) := \sum\limits_{n = 1}^\infty  {\frac{1}{{{n^p}}}}, \Re(p)>1. \end{eqnarray}
Some results for sums of (alternating) harmonic numbers may be seen in the works of \cite{BBG1994,BBG1995,BZB2008,F2005,M2014,M2009,S2015,S2008,S2010,So2011,S2014,SC2012,S2011,Xu2016,X2016} and references therein. Some explicit, and closely related results may also be seen in the well
presented papers \cite{SX2017,S2016,XMZ2016}. For example, in \cite{XMZ2016}, we give explicit formulas for the following type of Euler sums function
\[\sum\limits_{n = 1}^\infty  {\frac{{H_n^2}}{{{n^p}\left( {\begin{array}{*{20}{c}}
   {n + k}  \\
   k  \\
\end{array}} \right)}}{x^n}} ,\sum\limits_{n = 1}^\infty  {\frac{{H_n^{\left( m \right)}}}{{{n^p}\left( {\begin{array}{*{20}{c}}
   {n + k}  \\
   k  \\
\end{array}} \right)}}{x^n}},x\in [-1,0)\cup (0,1)\]
by using the method of partial fraction decomposition and integral representations of series. The purpose of the present paper is to establish closed form of harmonic number sums (1.3). Next, we begin with some basic notation. For $k\in \N$, ${\bf s}=(s_1,\ldots,s_k)\in (\mathbb{Z^*})^k\ (\mathbb{Z^*}:=\mathbb{Z}\setminus\{0\}=\{\pm1,\pm2,\ldots\})$, and a non-negative integer $n$, the multiple harmonic star sum is defined by (\cite{XMZ2016})
\[\zeta _n^ \star \left( {\bf{s}} \right) \equiv \zeta _n^ \star \left( {{s_1}, \ldots ,{s_k}} \right): = \sum\limits_{n \ge {n_1} \ge  \cdots  \ge {n_k} \ge 1} {\prod\limits_{j = 1}^k {n_j^{ - \left| {{s_j}} \right|}{\rm{sgn}}{{({s_j})}^{{n_j}}}} }.\tag{1.6} \]
Throughout the paper we will use ${\bar n}$ to denote a negative entry $s_j=-n$. For example,
\[\zeta _n^ \star \left( {\bar 2} \right) = \zeta _n^ \star \left( { - 2} \right),\zeta _n^ \star \left( {3,\bar 1} \right) = \zeta _n^ \star \left( {3, - 1} \right),\zeta _n^ \star \left( {\bar 2,2,\bar 1} \right) = \zeta _n^ \star \left( { - 2,2, - 1} \right).\]
We call $l({\bf s}):=k$ the depth of (1.6) and $\left| {\bf s} \right|: = \sum\limits_{j = 1}^k {\left| {{s_j}} \right|} $ the weight. For convenience we set ${\zeta_n ^ \star }\left( \emptyset  \right) = 1$ and ${\left\{ {{s_1}, \ldots ,{s_j}} \right\}_d}$ the set formed by repeating the composition $\left( {{s_1}, \ldots ,{s_j}} \right)$ $d$ times. When taking the limit $n\rightarrow \infty$ we get the so-called the star Euler sum
\[{\zeta ^ \star }\left( {\bf{s}} \right) = \mathop {\lim }\limits_{n \to \infty } \zeta _n^ \star \left( {\bf{s}} \right).\]
When ${\bf s}\in \N^k$ they are called the multiple zeta star value. It is obvious that
\[\zeta _n^ \star \left( m \right) = H_n^{\left( m \right)},m\in \N.\]
In this paper, we will prove that the alternating quadratic and cubic harmonic number sums ${{\overline W}_k}\left( {1,2;p} \right)$ and ${{\overline W}_k}\left( {{{\left\{ 1 \right\}}_3};p} \right)$ for $p=0,1$ can be expressed as a rational linear combination of products of single zeta values and multiple harmonic star sum of weight$\leq 4$ and depth $\leq 4$.
The main results of this paper as follow.
\begin{thm} For positive integer $k$, then the following identities hold:
\begin{align*}
\sum\limits_{n = 1}^\infty  {\frac{{{H_n}{H^{(2)} _n}}}{{n + k}}{{\left( { - 1} \right)}^{n + k}}}  =&  - \frac{5}{{16}}\zeta \left( 4 \right) - \frac{1}{4}\zeta \left( 2 \right){\ln ^2}2 + \frac{7}{8}\zeta \left( 3 \right)\ln 2 + \frac{7}{8}\zeta \left( 3 \right)\zeta _{k - 1}^ \star \left( {\bar 1} \right)\\
&- \frac{1}{4}\zeta \left( 3 \right){H_{k - 1}} - \frac{1}{2}\zeta \left( 2 \right)\zeta _{k - 1}^ \star \left( {\bar 2} \right) - \zeta _{k - 1}^ \star \left( {3,\bar 1} \right) + \zeta _{k - 1}^ \star \left( {1,2,\bar 1} \right)\\
&+ \frac{1}{2}\zeta \left( 2 \right)\zeta _{k - 1}^ \star \left( {1,\bar 1} \right) + \zeta _{k - 1}^ \star \left( {2,1,\bar 1} \right) + \ln 2\left\{ {\zeta _{k - 1}^ \star \left( {\bar 3} \right) - \zeta _{k - 1}^ \star \left( 3 \right)} \right\}\\
& - \frac{1}{2}\ln 2\zeta \left( 2 \right)\left\{ {\zeta _{k - 1}^ \star \left( {\bar 1} \right) - \zeta _{k - 1}^ \star \left( 1 \right)} \right\} + \frac{1}{2}{\ln ^2}2\left\{ {\zeta _{k - 1}^ \star \left( {\bar 2} \right) - \zeta _{k - 1}^ \star \left( 2 \right)} \right\}\\
& - \ln 2\left\{ {\zeta _{k - 1}^ \star \left( {2,\bar 1} \right) - \zeta _{k - 1}^ \star \left( {2,1} \right)} \right\} - \ln 2\left\{ {\zeta _{k - 1}^ \star \left( {1,\bar 2} \right) - \zeta _{k - 1}^ \star \left( {1,2} \right)} \right\},\tag{1.7}
\end{align*}
\begin{align*}
\sum\limits_{n = 1}^\infty  {\frac{{H_n^3}}{{n + k}}{{\left( { - 1} \right)}^{n + k}}}  =&  - \frac{5}{{16}}\zeta \left( 4 \right) + \frac{9}{8}\zeta \left( 3 \right)\ln 2 - \frac{3}{4}\zeta \left( 2 \right){\ln ^2}2 + \frac{1}{4}{\ln ^4}2 + \frac{9}{8}\zeta \left( 3 \right)\zeta _{k - 1}^ \star \left( {\bar 1} \right)\\
& - \frac{1}{2}\zeta \left( 2 \right)\zeta _{k - 1}^ \star \left( {\bar 2} \right) - \zeta _{k - 1}^ \star \left( {3,\bar 1} \right) + \frac{3}{2}\zeta \left( 2 \right)\zeta _{k - 1}^ \star \left( {1,\bar 1} \right) - \frac{3}{4}\zeta \left( 3 \right){H_{k - 1}}\\
& + 3\zeta _{k - 1}^ \star \left( {1,2,\bar 1} \right) + 3\zeta _{k - 1}^ \star \left( {2,1,\bar 1} \right) - 6\zeta _{k - 1}^ \star \left( {{{\left\{ 1 \right\}}_3},\bar 1} \right)\\
& + \ln 2\left\{ {\zeta _{k - 1}^ \star \left( {\bar 3} \right) - \zeta _{k - 1}^ \star \left( 3 \right)} \right\} - 3\ln 2\left\{ {\zeta _{k - 1}^ \star \left( {2,\bar 1} \right) - \zeta _{k - 1}^ \star \left( {2,1} \right)} \right\}\\
&- 3\ln 2\left\{ {\zeta _{k - 1}^ \star \left( {1,\bar 2} \right) - \zeta _{k - 1}^ \star \left( {1,2} \right)} \right\} + 6\ln 2\left\{ {\zeta _{k - 1}^ \star \left( {1,1,\bar 1} \right) - \zeta _{k - 1}^ \star \left( {{{\left\{ 1 \right\}}_3}} \right)} \right\}\\
&- 3{\ln ^2}2\left\{ {\zeta _{k - 1}^ \star \left( {1,\bar 1} \right) - \zeta _{k - 1}^ \star \left( {1,1} \right)} \right\} + \frac{3}{2}{\ln ^2}2\left\{ {\zeta _{k - 1}^ \star \left( {\bar 2} \right) - \zeta _{k - 1}^ \star \left( 2 \right)} \right\}\\
& + \left\{ {{{\ln }^3}2 - \frac{3}{2}\ln 2\zeta \left( 2 \right)} \right\}\left\{ {\zeta _{k - 1}^ \star \left( {\bar 1} \right) - \zeta _{k - 1}^ \star \left( 1 \right)} \right\}.\tag{1.8}
\end{align*}
\end{thm}
We will prove Theorem 1.1 in section 3.
\begin{thm} For integer $k\in\N$, we have
\begin{align*}
&{{\overline W}_k}\left( {1,2;0} \right) = \sum\limits_{r = 1}^k {{{\left( { - 1} \right)}^{r + 1}}r\left( {\begin{array}{*{20}{c}}
   k  \\
   r  \\
\end{array}} \right)\sum\limits_{n = 1}^\infty  {\frac{{{H_n}H_n^{\left( 2 \right)}}}{{n + r}}{{\left( { - 1} \right)}^{n + 1}}} }, \\
&{{\overline W}_k}\left( {{{\left\{ 1 \right\}}_3};0} \right) = \sum\limits_{r = 1}^k {{{\left( { - 1} \right)}^{r + 1}}r\left( {\begin{array}{*{20}{c}}
   k  \\
   r  \\
\end{array}} \right)\sum\limits_{n = 1}^\infty  {\frac{{H_n^3}}{{n + r}}{{\left( { - 1} \right)}^{n + 1}}} }, \\
&{{\overline W}_k}\left( {1,2;1} \right) = \sum\limits_{r = 1}^k {{{\left( { - 1} \right)}^{r + 1}}\left( {\begin{array}{*{20}{c}}
   k  \\
   r  \\
\end{array}} \right)\left\{ {\sum\limits_{n = 1}^\infty  {\frac{{{H_n}H_n^{\left( 2 \right)}}}{n}{{\left( { - 1} \right)}^{n + 1}} - \sum\limits_{n = 1}^\infty  {\frac{{{H_n}H_n^{\left( 2 \right)}}}{{n + r}}{{\left( { - 1} \right)}^{n + 1}}} } } \right\}}, \\
&{{\overline W}_k}\left( {{{\left\{ 1 \right\}}_3};1} \right) = \sum\limits_{r = 1}^k {{{\left( { - 1} \right)}^{r + 1}}\left( {\begin{array}{*{20}{c}}
   k  \\
   r  \\
\end{array}} \right)\left\{ {\sum\limits_{n = 1}^\infty  {\frac{{H_n^3}}{n}{{\left( { - 1} \right)}^{n + 1}} - \sum\limits_{n = 1}^\infty  {\frac{{H_n^3}}{{n + r}}{{\left( { - 1} \right)}^{n + 1}}} } } \right\}} .
\end{align*}
\end{thm}
\pf We consider the expansion
\[\frac{1}{{\prod\limits_{i = 1}^k {\left( {n + {a_i}} \right)} }} = \sum\limits_{j = 1}^k {\frac{{{A_j}}}{{n + {a_j}}}}\ \ (k\in \N_0; a_i \in \mathbb{C}\setminus\Z^-)\tag{1.9} \]
where
\[{A_j} = \mathop {\lim }\limits_{n \to  - {a_j}} \frac{{n + {a_j}}}{{\prod\limits_{i = 1}^k {\left( {n + {a_i}} \right)} }} = \prod\limits_{i = 1,i \ne j}^k {{{\left( {{a_i} - {a_j}} \right)}^{ - 1}}}.\tag{1.10}\]
Taking $a_i=a+i$ in (1.10) we obtain
\[{A_r} = {\left( { - 1} \right)^{r + 1}}\frac{r}
{{k!}}\left( {\begin{array}{*{20}{c}}
   k  \\
   r \\

 \end{array} } \right),\tag{1.11}\]
\[\frac{1}
{{\prod\limits_{i = 1}^k {\left( {n + a + i} \right)} }} = \sum\limits_{r = 1}^k {{{\left( { - 1} \right)}^{r + 1}}\frac{r}
{{k!}}\left( {\begin{array}{*{20}{c}}
   k  \\
   r  \\

 \end{array} } \right)\frac{1}
{{n + a + r}}} .\tag{1.12}\]
Furthermore, by using the equation (1.12) and the definition of binomial coefficient, letting $a=0$, we have the following expansions
\begin{align*}
&\frac{1}
{{\left( {\begin{array}{*{20}{c}}
   {n + k }  \\
   k  \\
\end{array} } \right)}} = \sum\limits_{r = 1}^k {{{\left( { - 1} \right)}^{r + 1}}r\left( {\begin{array}{*{20}{c}}
   k  \\
   r  \\
\end{array} } \right)\frac{1}
{{n + r}}} \ \ (k\in \N_0),\tag{1.13}
\end{align*}
Hence, by a direct calculation we may easily deduce the desired result. This completes the proof of Theorem 1.2. \hfill$\square$
\section{Some lemmas and theorems}
The following lemma will be useful in the development of the main theorem 1.1.
\begin{lem}(\cite{XMZ2016})
For integers $m,k\in \N$ and $x\in[-1,1)$, we have
\begin{align*}
&{\left( { - 1} \right)^m}m!\sum\limits_{n \ge m} {\frac{{s\left( {n + 1,m + 1} \right)}}{{\left( {n + k} \right)n!}}{x^{n + k}}}+ \frac{1}{{m + 1}}{\ln ^{m + 1}}\left( {1 - x} \right)\\
&=\sum\limits_{j = 1}^{m - 1} {{{\left( { - 1} \right)}^{j - 1}}j!\left( {\begin{array}{*{20}{c}}
   m  \\
   j  \\
\end{array}} \right)} {\ln ^{m - j}}\left( {1 - x} \right)\left( {\zeta _{k - 1}^ \star \left( {{{\{ 1\} }_{j + 1}};{{\{ 1\} }_j},x} \right) - \zeta _{k - 1}^ \star \left( {{{\{ 1\} }_{j + 1}}} \right)} \right)\\
& \quad- {\ln ^m}\left( {1 - x} \right)\left( {\zeta _{k - 1}^ \star \left( {1;x} \right) - \zeta _{k - 1}^ \star \left( 1 \right)} \right) - {\left( { - 1} \right)^m}m!\zeta _{k - 1}^ \star \left( {{{\{ 1\} }_{m + 1}};{{\{ 1\} }_m},x} \right).\tag{3.1}
\end{align*}
where ${s\left( {n,k} \right)}$ denotes the (unsigned) Stirling number of the first kind (\cite{L1974}),
\begin{align*}
& s\left( {n,1} \right) = \left( {n - 1} \right)!,s\left( {n,2} \right) = \left( {n - 1} \right)!{H_{n - 1}},s\left( {n,3} \right) = \frac{{\left( {n - 1} \right)!}}{2}\left[ {H_{n - 1}^2 - {H^{(2)} _{n - 1}}} \right],\\
&s\left( {n,4} \right) = \frac{{\left( {n - 1} \right)!}}{6}\left[ {H_{n - 1}^3 - 3{H_{n - 1}}{H^{(2)} _{n - 1}} + 2{H^{(3)} _{n - 1}}} \right], \\
&s\left( {n,5} \right) = \frac{{\left( {n - 1} \right)!}}{{24}}\left[ {H_{n - 1}^4 - 6{H^{(4)} _{n - 1}} - 6H_{n - 1}^2{H^{(2)} _{n - 1}} + 3(H^{(2)} _{n - 1})^2 + 8H_{n - 1}^{}{H^{(3)} _{n - 1}}} \right].
\end{align*}
The Stirling numbers ${s\left( {n,k} \right)}$ of the first kind satisfy a recurrence relation in the form
\[s\left( {n,k} \right) = s\left( {n - 1,k - 1} \right) + \left( {n - 1} \right)s\left( {n - 1,k} \right),\;\;n,k \in \N,\]
with $s\left( {n,k} \right) = 0,n < k,s\left( {n,0} \right) = s\left( {0,k} \right) = 0,s\left( {0,0} \right) = 1$.
The generating function of ${s\left( {n,k} \right)}$ is
\[{\ln ^k}\left( {1 - x} \right) = {\left( { - 1} \right)^k}k!\sum\limits_{n = k}^\infty  {\frac{{s\left( {n,k} \right)}}{{n!}}{x^n}} ,\: - 1 \le x < 1.\tag{3.2}\]
For $s_j>0$,  ${\zeta^\star _n}\left( {{s_1},{s_2}, \cdots ,{s_m};{x_1},{x_2}, \cdots ,{x_m}} \right)$ denotes the partial sums of multiple polylogarithm-star function defined by (\cite{XMZ2016})
\[{\zeta^\star _n}\left( {{s_1},{s_2}, \cdots ,{s_m};{x_1},{x_2}, \cdots ,{x_m}} \right) = \sum\limits_{1 \le {k_m} \le  \cdots  \le {k_1} \le n} {\prod\limits_{j = 1}^m {\frac{{x_j^{{k_j}}}}{{k_j^{{s_j}}}}} }. \]
\end{lem}
\begin{lem}(\cite{XMZ2016})
For integers $m,k\in \N$ and $x\in[-1,1)$, we have
\begin{align*}
\sum\limits_{n = 1}^\infty  {\frac{{H_n^{\left( m \right)}}}{{n + k}}{x^{n + k}}}  =& \sum\limits_{n = 1}^\infty  {\frac{{H_n^{\left( m \right)}}}{n}{x^n}}  - {\rm{L}}{{\rm{i}}_{m + 1}}\left( x \right) - \sum\limits_{j = 1}^{m - 1} {{{\left( { - 1} \right)}^{j - 1}}{\rm{L}}{{\rm{i}}_{m + 1 - j}}\left( x \right)} \zeta _{k - 1}^ \star \left( {j;x} \right)\\
& - {\left( { - 1} \right)^m}\ln \left( {1 - x} \right)\left\{ {\zeta _{k - 1}^ \star \left( {m;x} \right) - \zeta _{k - 1}^ \star \left( m \right)} \right\} + {\left( { - 1} \right)^m}\zeta _{k - 1}^ \star \left( {m,1;1,x} \right),\tag{3.3}
\end{align*}
where ${\rm Li}{_p}\left( x \right)$ denotes the polylogarithm function defined for $\left| x \right| \le 1$ by the series
\[{\rm Li}{_p}\left( x \right) := \sum\limits_{n = 1}^\infty  {\frac{{{x^n}}}{{{n^p}}}}, \Re(p)>1 .\tag{3.4}\]
\end{lem}
\begin{lem}(\cite{SX2017,X2016})
For integers $n\geq1, k\geq0$, then
\[\int\limits_0^1 {{t^{n - 1}}{{\ln }^k}\left( {1 - t} \right)} dt = {\left( { - 1} \right)^k}\frac{{{Y_k}\left( n \right)}}{n},{Y_0}\left( n \right) = 1,\tag{3.5}\]
where ${Y_k}\left( n \right) = {Y_k}\left( {{H _n},1!{H^{(2)} _n},2!{H^{(3)} _n}, \cdots ,\left( {r - 1} \right)!{H^{(r)} _n}, \cdots } \right)$, ${Y_k}\left( {{x_1},{x_2}, \cdots } \right)$ stands for the complete exponential Bell polynomial defined by (see \cite{L1974})
\[\exp \left( {\sum\limits_{m \ge 1}^{} {{x_m}\frac{{{t^m}}}{{m!}}} } \right) = 1 + \sum\limits_{k \ge 1}^{} {{Y_k}\left( {{x_1},{x_2}, \cdots } \right)\frac{{{t^k}}}{{k!}}}.\tag{3.6}\]
\end{lem}
 From the definition of the complete exponential Bell polynomial, we have
$${Y_1}\left( n \right) = {H_n},{Y_2}\left( n \right) = H_n^2 + {H^{(2)} _n},{Y_3}\left( n \right) =  H_n^3+ 3{H_n}{H^{(2)} _n}+ 2{H^{(3)} _n},$$
\[{Y_4}\left( n \right) = H_n^4 + 8{H_n}{H^{(3)} _n} + 6H_n^2{H^{(2)} _n} + 3(H^{(2)} _n)^2 + 6{H^{(4)} _n}.\]
\begin{lem}(\cite{SX2017})
For integer $m \geq 1$, and $-1<x<1$ , then
\[\sum\limits_{n = 1}^\infty  {{H_n}{H^{(m)} _n}{x^n}}  = \frac{1}{{1 - x}}\left\{ {\sum\limits_{n = 1}^\infty  {\frac{{{H_n}}}{{{n^m}}}{x^n}}  - \sum\limits_{n = 1}^\infty  {\frac{1}{{{n^m}}}\left( {\sum\limits_{k = 1}^n {\frac{{{x^k}}}{k}} } \right)}  - \zeta \left( m \right)\ln \left( {1 - x} \right)} \right\}.\tag{3.7}\]
\end{lem}
\begin{thm} For any real $x\in (-1,1)$, then the following identity holds:
\[\sum\limits_{n = 1}^\infty  {{H_n}H_n^{\left( 2 \right)}{x^n}}  = \frac{1}{{1 - x}}\left\{ {2{\rm Li}{_3}\left( x \right) - \ln \left( {1 - x} \right){\rm Li}{_2}\left( x \right) - \sum\limits_{n = 1}^\infty  {\frac{{{H_n}}}{{{n^2}}}{x^n}} } \right\}.\tag{3.8}\]
\end{thm}
\pf To prove identity (3.8), we consider the nested sum
\[\sum\limits_{n = 1}^\infty  {\frac{{{y^n}}}{{{n^m}}}\left( {\sum\limits_{k = 1}^n {\frac{{{x^k}}}{{{k^p}}}} } \right)} ,\;x,y \in \left[ { - 1,1} \right),\;m,p \in \N.\]
By taking the sum over complementary pairs of summation indices, we obtain a simple reflection formula
\[\sum\limits_{n = 1}^\infty  {\frac{{{y^n}}}{{{n^m}}}\left( {\sum\limits_{k = 1}^n {\frac{{{x^k}}}{{{k^p}}}} } \right)}  + \sum\limits_{n = 1}^\infty  {\frac{{{x^n}}}{{{n^p}}}\left( {\sum\limits_{k = 1}^n {\frac{{{y^k}}}{{{k^m}}}} } \right)}  = {\rm Li}{_p}\left( x \right){\rm Li}{_m}\left( y \right) + {\rm Li}{_{p + m}}\left( {xy} \right).\tag{3.9}\]
Setting $p=1,m=2,y=1$ in above equation we get
\[\sum\limits_{n = 1}^\infty  {\frac{1}{{{n^2}}}\left( {\sum\limits_{k = 1}^n {\frac{{{x^k}}}{k}} } \right)}  + \sum\limits_{n = 1}^\infty  {\frac{{H_n^{\left( 2 \right)}}}{n}{x^n}}  =  - \ln \left( {1 - x} \right)\zeta \left( 2 \right) + {\rm{L}}{{\rm{i}}_3}\left( x \right).\tag{3.10}\]
On the other hand, by the definition of polylogarithm function and Cauchy product of power series, we have
\begin{align*}
\sum\limits_{n = 1}^\infty  {\frac{{H_n^{\left( 2 \right)}}}{n}{x^n}} & =\int\limits_0^x {\frac{{{\rm{L}}{{\rm{i}}_2}\left( t \right)}}{{t\left( {1 - t} \right)}}dt}  = \int\limits_0^x {\frac{{{\rm{L}}{{\rm{i}}_2}\left( t \right)}}{t}dt}  + \int\limits_0^x {\frac{{{\rm{L}}{{\rm{i}}_2}\left( t \right)}}{{1 - t}}dt} \\
& = {\rm{L}}{{\rm{i}}_3}\left( x \right) - \ln \left( {1 - x} \right){\rm{L}}{{\rm{i}}_2}\left( x \right) - \int\limits_0^x {\frac{{{{\ln }^2}\left( {1 - t} \right)}}{t}} dt\\
& = 3{\rm{L}}{{\rm{i}}_3}\left( x \right) - \ln \left( {1 - x} \right){\rm{L}}{{\rm{i}}_2}\left( x \right) - 2\sum\limits_{n = 1}^\infty  {\frac{{{H_n}}}{{{n^2}}}{x^n}} .\tag{3.11}
\end{align*}
Then, substituting (3.11) into (3.10) yields
\[\sum\limits_{n = 1}^\infty  {\frac{1}{{{n^2}}}\left( {\sum\limits_{k = 1}^n {\frac{{{x^k}}}{k}} } \right)}  = 2\sum\limits_{n = 1}^\infty  {\frac{{{H_n}}}{{{n^2}}}{x^n}}  + \ln \left( {1 - x} \right){\rm{L}}{{\rm{i}}_2}\left( x \right) - 2{\rm{L}}{{\rm{i}}_3}\left( x \right) - \ln \left( {1 - x} \right)\zeta \left( 2 \right).\tag{3.12}\]
Hence, taking $m=2$ in Lemma 2.4 and combining formula (3.12) we may deduce the desired result.
 The proof of Theorem 2.5 is thus completed.\hfill$\square$
 \begin{thm}
If $m\geq1$ is a integer and $z \in \left[ {0,1} \right]$
, then we have
\begin{align*}\int\limits_0^z {\frac{{{{\ln }^m}\left( {1 + x} \right)}}{x}} dx = &\frac{1}{{m + 1}}{\ln ^{m + 1}}\left( {1 + z} \right) + m!\left( {\zeta \left( {m + 1} \right) - {\rm Li}{_{m + 1}}\left( {\frac{1}{{1 + z}}} \right)} \right) \\&- m!\sum\limits_{j = 1}^m {\frac{{{{\ln }^{m - j + 1}}\left( {1 + z} \right)}}{{\left( {m - j + 1} \right)!}}} {\rm Li}{_j}\left( {\frac{1}{{1 + z}}} \right).\tag{3.13}
\end{align*}
\end{thm}
\pf We note that the integral in (3.13) can be rewritten as
\begin{align*}
\int\limits_0^z {\frac{{{{\ln }^m}\left( {1 + x} \right)}}{x}} dx\mathop  = \limits^{t = 1 + x}  & \int\limits_1^{1 + z} {\frac{{{{\ln }^m}t}}{{t - 1}}} dt\mathop  = \limits^{u = {t^{ - 1}}} {\left( { - 1} \right)^{m + 1}}\int\limits_1^{{{(1 + z)^{-1}}}} {\frac{{{{\ln }^m}u}}{{u - {u^2}}}} du\\
           =& {\left( { - 1} \right)^{m + 1}}\left\{ {\int\limits_1^{{{(1 + z)^{-1}}}} {\frac{{{{\ln }^m}u}}{u}} du + \int\limits_1^{{{(1 + z)^{-1}}}} {\frac{{{{\ln }^m}u}}{{1 - u}}} du} \right\}\\
           =&\frac{1}{{m + 1}}{\ln ^{m + 1}}\left( {1 + z} \right) + {\left( { - 1} \right)^{m + 1}}\int\limits_1^{{{(1 + z)^{-1}}}} {\frac{{{{\ln }^m}u}}{{1 - u}}} du .\tag{3.14}
\end{align*}
We use the elementary integral identity
$$\int {{t^{n - 1}}{{\left( {\ln t} \right)}^m}dt}  = {t^n}\left\{ {\frac{{{{\ln }^m}t}}{n} - \sum\limits_{l = 1}^m {\frac{{{{\left( { - 1} \right)}^{l + 1}}{{\ln }^{m - l}}\left( {t} \right){{\left( m \right)}_l}}}{{{n^{l + 1}}}}} } \right\},\eqno(3.15)$$
Here $(m)_l=m(m-1)\cdots(m-l+1)$.
Then yields
\begin{align*}
\int\limits_1^{{{(1 + z)^{-1}}}} {\frac{{{{\ln }^m}u}}{{1 - u}}} du = \sum\limits_{n = 1}^\infty  {\int\limits_1^{{{(1 + z)^{-1}}}} {{u^{n - 1}}{{\ln }^m}u} du} = & {\left( { - 1} \right)^{m + 1}}m!\left( {\zeta \left( {m + 1} \right) - {\rm Li}{_{m + 1}}\left( {{{(1 + z)^{-1}}}} \right)} \right)\\
           & + m!{\left( { - 1} \right)^m}\sum\limits_{j = 1}^m {\frac{{{{\ln }^{m - j + 1}}\left( {1 + z} \right)}}{{\left( {m - j + 1} \right)!}}} {\rm Li}{_j}\left( {{{(1 + z)^{-1}}}} \right).\tag{3.16}
\end{align*}
Substituting (3.16) into (3.14) yields the desired result.\hfill$\square$  \\
Setting $p=3,4$ in (3.13), we get
\begin{align*}
 &\int\limits_0^1 {\frac{{{{\ln }^3}\left( {1 + x} \right)}}{x}} dx = 6\zeta \left( 4 \right) + \frac{3}{2}\zeta \left( 2 \right){\ln ^2}2 - \frac{1}{4}{\ln ^4}2 - \frac{{21}}{4}\zeta \left( 3 \right)\ln 2 - 6{\rm{L}}{{\rm{i}}_4}\left( {\frac{1}{2}} \right), \\
 &\int\limits_0^1 {\frac{{{{\ln }^4}\left( {1 + x} \right)}}{x}} dx =  - 24{\rm{L}}{{\rm{i}}_5}\left( {\frac{1}{2}} \right) - 24\ln 2{\rm{L}}{{\rm{i}}_4}\left( {\frac{1}{2}} \right) - \frac{4}{5}{\ln ^5}2 - \frac{{21}}{2}\zeta \left( 3 \right){\ln ^2}2 + 24\zeta \left( 5 \right) + 4\zeta \left( 2 \right){\ln ^3}2.
\end{align*}
\begin{thm} For integer $m\in \N$, then
\[\sum\limits_{n = 1}^\infty  {\frac{{{Y_m}\left( n \right)}}{n}{{\left( { - 1} \right)}^{n - 1}}}  = m!{\rm{L}}{{\rm{i}}_{m + 1}}\left( {\frac{1}{2}} \right).\tag{3.17}\]
\end{thm}
\pf First, by a direct calculation, we get
\[\int\limits_0^1 {\frac{{{{\ln }^m}\left( {1 - x} \right)}}{{1 + x}}} dx = {\left( { - 1} \right)^m}m!{\rm{L}}{{\rm{i}}_{m + 1}}\left( {\frac{1}{2}} \right).\tag{3.18} \]
Then by applying the known result
\[\frac{1}{{1 + x}} = \sum\limits_{n = 1}^\infty  {{{\left( { - 1} \right)}^{n - 1}}{x^{n-1}}} ,x \in \left( { - 1,1} \right),\]
in (3.18), we obtain
\begin{align*}
\int\limits_0^1 {\frac{{{{\ln }^m}\left( {1 - x} \right)}}{{1 + x}}} dx =& \sum\limits_{n = 1}^\infty  {{{\left( { - 1} \right)}^{n - 1}}} \int\limits_0^1 {{x^{n - 1}}{{\ln }^m}\left( {1 - x} \right)} dx\\
& = {\left( { - 1} \right)^m}\sum\limits_{n = 1}^\infty  {\frac{{{Y_m}\left( n \right)}}{n}{{\left( { - 1} \right)}^{n - 1}}} .\tag{3.19}
\end{align*}
Therefore, combining (3.18) with (3.19), we obtain the formula (3.17). This completes the proof of Theorem 2.7.\hfill$\square$
\begin{cor} The following identities hold:
\begin{align*}
&\sum\limits_{n = 1}^\infty  {\frac{{H_n^3}}{n}{{\left( { - 1} \right)}^{n - 1}}}  = \frac{5}{8}\zeta \left( 4 \right) + \frac{3}{4}\zeta \left( 2 \right){\ln ^2}2 - \frac{1}{4}{\ln ^4}2 - \frac{9}{8}\zeta \left( 3 \right)\ln 2,\tag{3.20}\\
&\sum\limits_{n = 1}^\infty  {\frac{{H_n^{}H^{(2)}_n}}{n}} {\left( { - 1} \right)^{n - 1}} = 2{\rm{L}}{{\rm{i}}_4}\left( {\frac{1}{2}} \right) + \frac{1}{{12}}{\ln ^4}2{\rm{ + }}\frac{7}{8}\zeta \left( 3 \right)\ln 2 - \frac{1}{4}\zeta \left( 2 \right){\ln ^2}2 - \zeta \left( 4 \right).\tag{3.21}
\end{align*}
\end{cor}
\pf Taking $m=3$ in Theorem 2.6 and Theorem 2.7, we have
\begin{align*}
&\sum\limits_{n = 1}^\infty  {\frac{{H_n^3 - 3{H_n}{H^{(2)} _n} + 2{H^{(3)} _n}}}{n}{{\left( { - 1} \right)}^{n - 1}} = \int\limits_0^1 {\frac{{{{\ln }^3}\left( {1 + x} \right)}}{x}dx} }  - \frac{1}{4}{\ln ^4}2,\tag{3.22}\\
&\sum\limits_{n = 1}^\infty  {\frac{{H_n^3 + 3{H_n}{H^{(2)} _n} + 2{H^{(3)} _n}}}{n}{{\left( { - 1} \right)}^{n - 1}}}  = 6{\rm{L}}{{\rm{i}}_4}\left( {\frac{1}{2}} \right).\tag{3.23}
\end{align*}
From \cite{FS1998} we deduce that
\[\sum\limits_{n = 1}^\infty  {\frac{1}{{{n^3}}}\left( {\sum\limits_{k = 1}^n {\frac{{{{\left( { - 1} \right)}^{k - 1}}}}{k}} } \right)}  = \frac{7}{4}\zeta \left( 3 \right)\ln 2 - \frac{5}{{16}}\zeta \left( 4 \right).\tag{3.24}\]
Then substituting (3.24) into (3.9) with $m=3,p=1,x=-1,y=1$ yields
\[\sum\limits_{n = 1}^\infty  {\frac{{H_n^{\left( 3 \right)}}}{n}{{\left( { - 1} \right)}^{n - 1}}}  = \frac{{19}}{{16}}\zeta (4) - \frac{3}{4}\zeta \left( 3 \right)\ln 2.\tag{3.25}\]
Thus, combining formulas (3.22), (3.23) and (3.25) we obtain the results.\hfill$\square$
\begin{cor} For positive integer $k$, then we have
\begin{align*}\sum\limits_{n = 1}^\infty  {\frac{{H_n^3 - 3{H_n}{H^{(2)} _n} + 2{H^{(3)} _n}}}{{n + k}}{{\left( { - 1} \right)}^{n + k}}}  =& \frac{1}{4}{\ln ^4}2 - 6\zeta _{k - 1}^ \star \left( {{{\left\{ 1 \right\}}_3},\bar 1} \right) + {\ln ^3}2\left\{ {\zeta _{k - 1}^ \star \left( {\bar 1} \right) - \zeta _{k - 1}^ \star \left( 1 \right)} \right\}\\
& - 3{\ln ^2}2\left\{ {\zeta _{k - 1}^ \star \left( {1,\bar 1} \right) - \zeta _{k - 1}^ \star \left( {1,1} \right)} \right\}\\& + 6\ln 2\left\{ {\zeta _{k - 1}^ \star \left( {1,1,\bar 1} \right) - \zeta _{k - 1}^ \star \left( {{{\left\{ 1 \right\}}_3}} \right)} \right\},\tag{3.26}
\end{align*}
\begin{align*}
\sum\limits_{n = 1}^\infty  {\frac{{H_n^{\left( 3 \right)}}}{{n + k}}{{\left( { - 1} \right)}^{n + k}}}  =&  - \frac{5}{{16}}\zeta \left( 4 \right) + \frac{3}{4}\zeta \left( 3 \right)\ln 2 + \ln 2\left\{ {\zeta _{k - 1}^ \star \left( {\bar 3} \right) - \zeta _{k - 1}^ \star \left( 3 \right)} \right\} + \frac{3}{4}\zeta \left( 3 \right)\zeta _{k - 1}^ \star \left( {\bar 1} \right)\\
& - \frac{1}{2}\zeta \left( 2 \right)\zeta _{k - 1}^ \star \left( {\bar 2} \right) - \zeta _{k - 1}^ \star \left( {3,\bar 1} \right).\tag{3.27}
\end{align*}
\end{cor}
\pf Taking $m=3,x=-1$ in Lemma 2.1 and Lemma 2.2 we obtain the result.\hfill$\square$
\section{Proof of Theorem 1.1}
In (3.8), changing $x$ to $t$, then multiplying it by $t^{k-1}-t^{-1}$ and integrating over $(0,x)$, we have
\begin{align*}
&\sum\limits_{n = 1}^\infty  {\left\{ {\frac{{{H_n}H_n^{\left( 2 \right)}}}{{n + k}}{x^{n + k}} - \frac{{{H_n}H_n^{\left( 2 \right)}}}{n}{x^n}} \right\}} \\
&=  - \sum\limits_{i = 1}^k {\int\limits_0^x {{t^{i - 2}}\left\{ {2{\rm{L}}{{\rm{i}}_3}\left( t \right) - \ln \left( {1 - t} \right){\rm{L}}{{\rm{i}}_2}\left( t \right) - \sum\limits_{n = 1}^\infty  {\frac{{{H_n}}}{{{n^2}}}{t^n}} } \right\}dt} } \\
& =  - \left\{ {2\int\limits_0^x {\frac{{{\rm{L}}{{\rm{i}}_3}\left( t \right)}}{t}dt - \int\limits_0^x {\frac{{\ln \left( {1 - t} \right){\rm{L}}{{\rm{i}}_2}\left( t \right)}}{t}dt}  - \sum\limits_{n = 1}^\infty  {\frac{{{H_n}}}{{{n^2}}}\int\limits_0^x {{t^{n - 1}}dt} } } } \right\}\\
&\quad - \sum\limits_{i = 1}^{k - 1} {\int\limits_0^x {{t^{i - 1}}\left\{ {2{\rm{L}}{{\rm{i}}_3}\left( t \right) - \ln \left( {1 - t} \right){\rm{L}}{{\rm{i}}_2}\left( t \right) - \sum\limits_{n = 1}^\infty  {\frac{{{H_n}}}{{{n^2}}}{t^n}} } \right\}dt} } \\
&=  - 2{\rm{L}}{{\rm{i}}_4}\left( x \right) - \frac{1}{2}{\rm{Li}}_2^2\left( x \right) + \sum\limits_{n = 1}^\infty  {\frac{{{H_n}}}{{{n^3}}}{x^n}} \\
&\quad + \sum\limits_{i = 1}^{k - 1} {\int\limits_0^x {{t^{i - 1}}\ln \left( {1 - t} \right){\rm{L}}{{\rm{i}}_2}\left( t \right)dt} }  + \sum\limits_{i = 1}^{k - 1} {\sum\limits_{n = 1}^\infty  {\frac{{{H_n}}}{{{n^2}\left( {n + i} \right)}}{x^{n + i}}} } \\
&\quad - 2\sum\limits_{i = 1}^{k - 1} {\left\{ {\frac{{{x^i}}}{i}{\rm{L}}{{\rm{i}}_3}\left( x \right) - \frac{{{x^i}}}{{{i^2}}}{\rm{L}}{{\rm{i}}_2}\left( x \right) - \frac{1}{{{i^3}}}\ln \left( {1 - x} \right)\left( {{x^i} - 1} \right) + \frac{1}{{{i^3}}}\sum\limits_{j = 1}^i {\frac{{{x^j}}}{j}} } \right\}} .\tag{4.1}
\end{align*}
In \cite{XMZ2016}, we deduce the following identity
\[\int_0^x {{t^{n - 1}}\ln \left( {1 - t} \right)} dt = \frac{1}{n}\left\{ {{x^n}\ln \left( {1 - x} \right) - \sum\limits_{j = 1}^n {\frac{{{x^j}}}{j}}  - \ln \left( {1 - x} \right)} \right\}.\tag{4.2}\]
Hence, by using (4.2) with the help of Theorem 2.2 in the reference \cite{XMZ2016}, then
\begin{align*}
&\int\limits_0^x {{t^{i - 1}}\ln \left( {1 - t} \right){\rm{L}}{{\rm{i}}_2}\left( t \right)dt}  = \int\limits_0^x {{\rm{L}}{{\rm{i}}_2}\left( t \right)d\left\{ {\frac{1}{i}\left[ {\left( {{t^i} - 1} \right)\ln \left( {1 - t} \right) - \sum\limits_{j = 1}^i {\frac{{{t^j}}}{j}} } \right]} \right\}} \\
& = \frac{1}{i}\left\{ {\left( {{x^i} - 1} \right)\ln \left( {1 - x} \right) - \sum\limits_{j = 1}^i {\frac{{{x^j}}}{j}} } \right\}{\rm{L}}{{\rm{i}}_2}\left( x \right) + \frac{1}{i}\int\limits_0^x {{t^{i - 1}}{{\ln }^2}\left( {1 - t} \right)dt} \\
&\quad - \frac{1}{i}\int\limits_0^x {\frac{{{{\ln }^2}\left( {1 - t} \right)}}{t}dt}  - \frac{1}{i}\sum\limits_{j = 1}^i {\frac{1}{j}\int\limits_0^x {{t^{j - 1}}\ln \left( {1 - t} \right)dt} } \\
& = \frac{1}{i}\left\{ {\left( {{x^i} - 1} \right)\ln \left( {1 - x} \right) - \zeta _i^ \star \left( {1;x} \right)} \right\}{\rm{L}}{{\rm{i}}_2}\left( x \right) + \frac{1}{{{i^2}}}{\ln ^2}\left( {1 - x} \right)\left( {{x^i} - 1} \right)\\
&\quad - \frac{2}{{{i^2}}}\ln \left( {1 - x} \right)\left\{ {\zeta _i^ \star \left( {1;x} \right) - \zeta _i^ \star \left( 1 \right)} \right\} + \frac{2}{{{i^2}}}\zeta _i^ \star \left( {1,1;1,x} \right)\\
&\quad - \frac{1}{i}\ln \left( {1 - x} \right)\left\{ {\left\{ {\zeta _i^ \star \left( {2;x} \right) - \zeta _i^ \star \left( 2 \right)} \right\}} \right\} + \frac{1}{i}\zeta _i^ \star \left( {2,1;1,x} \right)\\
&\quad - \frac{2}{i}\left\{ {\sum\limits_{n = 1}^\infty  {\frac{{{H_n}}}{{{n^2}}}{x^n} - {\rm{L}}{{\rm{i}}_3}\left( x \right)} } \right\}.\tag{4.3}
\end{align*}
Moreover, by a direct calculation we obtain
\begin{align*}
&\sum\limits_{i = 1}^{k - 1} {\sum\limits_{n = 1}^\infty  {\frac{{{H_n}}}{{{n^2}\left( {n + i} \right)}}{{\left( { - 1} \right)}^{n + i}} = } } \sum\limits_{i = 1}^{k - 1} {{{\left( { - 1} \right)}^{i - 1}}\sum\limits_{n = 1}^\infty  {{H_n}{{\left( { - 1} \right)}^{n - 1}}\left\{ {\frac{1}{i} \cdot \frac{1}{{{n^2}}} - \frac{1}{{{i^2}}} \cdot \frac{1}{n} + \frac{1}{{{i^2}}} \cdot \frac{1}{{n + i}}} \right\}} } \\
&= \sum\limits_{i = 1}^{k - 1} {\frac{{{{\left( { - 1} \right)}^{i - 1}}}}{i}} \sum\limits_{n = 1}^\infty  {\frac{{{H_n}}}{{{n^2}}}{{\left( { - 1} \right)}^{n - 1}}}  - \sum\limits_{i = 1}^{k - 1} {\frac{{{{\left( { - 1} \right)}^{i - 1}}}}{{{i^2}}}} \sum\limits_{n = 1}^\infty  {\frac{{{H_n}}}{n}{{\left( { - 1} \right)}^{n - 1}}}  + \sum\limits_{i = 1}^{k - 1} {\frac{{{{\left( { - 1} \right)}^{i - 1}}}}{{{i^2}}}} \sum\limits_{n = 1}^\infty  {\frac{{{H_n}}}{{n + i}}{{\left( { - 1} \right)}^{n - 1}}} \\
& =  - \frac{5}{8}\zeta \left( 3 \right)\zeta _{k - 1}^ \star \left( {\bar 1} \right) + \frac{{\zeta \left( 2 \right) - {{\ln }^2}2}}{2}\zeta _{k - 1}^ \star \left( {\bar 2} \right) + \ln 2\left\{ {\zeta _{k - 1}^ \star \left( {2,\bar 1} \right) - \zeta _{k - 1}^ \star \left( {2,1} \right)} \right\}\\
&\quad- \ln 2\left\{ {\zeta _{k - 1}^ \star \left( {\bar 3} \right) - \zeta _{k - 1}^ \star \left( 3 \right)} \right\} + \frac{1}{2}{\ln ^2}2\zeta _{k - 1}^ \star \left( {\bar 2} \right) - \zeta _{k - 1}^ \star \left( {2,1,\bar 1} \right) + \zeta _{k - 1}^ \star \left( {3,\bar 1} \right).\tag{4.4}
\end{align*}
Hence, letting $x=-1$ in (4.1), (4.3) and combining identities (3.21), (3.26), (3.27) and (4.4),
and using the following formula (\cite{FS1998})
\[\sum\limits_{n = 1}^\infty  {\frac{{{H_n}{}}}{{{n^3}}}}{\left( { - 1} \right)}^{n - 1}  =  - 2{\rm Li}{_4}\left( {\frac{1}{2}} \right) + \frac{{11{}}}{{4}} \zeta(4) + \frac{{{1}}}{{2}}\zeta(2){\ln ^2}2 - \frac{1}{{12}}{\ln ^4}2 - \frac{7}{4}\zeta \left( 3 \right)\ln 2,\]
by a simple calculation, we can prove the Theorem 1.1.\hfill$\square$\\
Taking $k=2$ in Theorem 1.1, we can get the following results.
\begin{cor} The following identities hold:
\begin{align*}
&\sum\limits_{n = 1}^\infty  {\frac{{H_n^3}}{{n + 2}}{{\left( { - 1} \right)}^n}}  =  - \frac{5}{{16}}\zeta \left( 4 \right) + \frac{9}{8}\zeta \left( 3 \right)\ln 2 - \frac{3}{4}\zeta \left( 2 \right){\ln ^2}2 + \frac{1}{4}{\ln ^4}2- 2\ln 2\\
&\quad\quad \quad\quad \quad\quad \quad\quad  + 3{\ln ^2}2 - 2{\ln ^3}2 + 3\zeta \left( 2 \right)\ln 2 - \frac{{15}}{8}\zeta \left( 3 \right) - \zeta \left( 2 \right) + 1,\\
&\sum\limits_{n = 1}^\infty  {\frac{{{H_n}H_n^{\left( 2 \right)}}}{{n + 2}}{{\left( { - 1} \right)}^n}}  =  - \frac{5}{{16}}\zeta \left( 4 \right) - \frac{1}{4}\zeta \left( 2 \right){\ln ^2}2 + \frac{7}{8}\zeta \left( 3 \right)\ln 2\\
&\quad\quad \quad \quad \quad \quad \quad \quad \quad  - \frac{9}{8}\zeta \left( 3 \right) + \zeta \left( 2 \right)\ln 2 + 2\ln 2 - {\ln ^2}2 - 1.
\end{align*}
\end{cor}
\section{Conclusion}
From Theorem 1.1, Theorem 1.2 and Corollary 2.8, we obtain the following description.
\begin{thm} For $p=0,1$, then the alternating quadratic and cubic harmonic number sums ${{\overline W}_k}\left( {1,2;p} \right)$ and ${{\overline W}_k}\left( {{{\left\{ 1 \right\}}_3};p} \right)$ can be expressed as a rational linear combination of products of single zeta values and multiple harmonic star sum of weight$\leq 4$ and depth $\leq 4$.
\end{thm}
{\bf Acknowledgments.} The authors would like to thank the anonymous
referee for his/her helpful comments, which improve the presentation
of the paper.

 \end{document}